\documentclass[12pt]{article}
\usepackage[pctex32]{graphics}
\usepackage{amsthm,amsmath,amssymb,amscd,verbatim,epsfig}
\usepackage{amssymb}
\usepackage{graphicx}
\usepackage{epsfig}
\usepackage{graphics,pifont}
\newcommand{\beq}{\begin{equation}}
\newcommand{\eeq}{\end{equation}}

\date{}

\newcommand{\f}{\frac}

\begin{document}

\title{Congruence\ Identities\ Arising\ From\ Dynamical\ Systems}
\author{Bau-Sen Du \\ [.5cm]
Institute of Mathematics \\
Academia Sinica \\
Taipei 11529, Taiwan \\
dubs@math.sinica.edu.tw \\
(Appl. Math. Letters, 12(1999), 115-119) \\}
\maketitle


\begin{abstract}
By counting the numbers of periodic points of all periods for some interval maps, we obtain infinitely many new congruence identities in number theory.
\end{abstract}

Let $S$ be a nonempty set and let $f$ be a map from $S$ into itself.  For every positive integer $n$, we define the $n^{th}$ 
iterate of $f$ by letting $f^1 = f$ and $f^n = f \circ f^{n-1}$ for $n \ge 2$.  For $y \in S$, we call the set $\{ \, f^k(y) \,: \, k \ge 0 \, \}$ the orbit of $y$ under $f$.  If $f^m(y) = y$ for some positive integer $m$, we call $y$ a periodic point of $f$ and call the smallest such positive integer $m$ the least period of $y$ under $f$.  We also call periodic points of least period 1 fixed points.  It is clear that if  $y$ is a periodic point of $f$ with least period $m$, then, for every integer $1 \le k \le m-1$, $f^k(y)$ is also a periodic point of $f$ with least period $m$ and they are all distinct.  So, every periodic orbit of $f$ with least period $m$ consists of exactly $m$ points.  Since distinct periodic orbits of $f$ are pairwise disjoint, the number (if finite) of distinct periodic points of $f$ with least period $m$ is divisible by $m$ and the quotient equals the number of distinct periodic orbits of $f$ with least period $m$.  Therefore, if there is a way to find the numbers of periodic points of all periods for a map, then we obtain infinitely many congruence identities in number theory.  This is an interesting application of dynamical systems theory to number theory which is not found in {\bf{\cite{la, sl}}}.

Let $\phi(m)$ be an integer-valued function defined on the set of all positive integers.  If $m=p_1^{k_1}p_2^{k_2} \cdots p_r^{k_r}$, where the $p_i$'s are distinct prime numbers, $r$ and $k_i$'s are positive integers, we let $\Phi_1(1, \phi) = \phi(1)$ and let $\Phi_1(m, \phi) =$
$$
\phi(m)-\sum_{i=1}^r \phi(\f m{p_i})+\sum_{i_1<i_2} \phi(\f m{p_{i_1}p_{i_2}})
- \sum_{i_1<i_2<i_3} \phi(\f m{p_{i_1}p_{i_2}p_{i_3}}) + \cdots 
+ (-1)^r \phi(\f m{p_1p_2 \cdots p_r}),
$$
\noindent
where the summation $\sum_{i_1<i_2< \cdots < i_j}$ is taken over all integers $i_1, i_2, \cdots, i_j$ with $1 \le i_1 < i_2 < \cdots < i_j \le r$.  If $m = 2^{k_0}p_1^{k_1}p_2^{k_2} \cdots p_r^{k_r}$, where the $p_i$'s are distinct odd prime numbers, and $k_0 \ge 0, r\ge 1$, and the $k_i$'s $\ge 1$ are integers, we let $\Phi_2(m, \phi) =$
$$
\phi(m)-\sum_{i=1}^r \phi(\f m{p_i})+\sum_{i_1<i_2} \phi(\f m{p_{i_1}p_{i_2}})
- \sum_{i_1<i_2<i_3} \phi(\f m{p_{i_1}p_{i_2}p_{i_3}}) + \cdots 
+ (-1)^r \phi(\f m{p_1p_2 \cdots p_r}), 
$$
\noindent
If $m = 2^k$, where $k \ge 0$ is an integer, we let $\Phi_2(m, \phi) = \phi(m) - 1$.  

Let $f$ be a map from the set $S$ into itself.  For every positive integer $m = p_1^{k_1}p_2^{k_2} \cdots p_r^{k_r}$, where $p_i$'s and $k_i$'s are defined as above, if $\phi (m)$ represents the number of distinct solutions of the equation $f^m(x) = x$ (i.e. the number of fixed points of $f^m(x)$) in $S$, then in the above formula for $\Phi_1(m, \phi)$, the periodic points of $f$ with least period $\f m{p_{i_1}^{t_{i_1}}p_{i_2}^{t_{i_2}} \cdots p_{i_j}^{t_{i_j}}} < m$, where $1 \le t_{i_s} \le k_{i_s}$, $1 \le s \le j$ are integers, have been counted 
\newline
\indent
\indent
\indent
\indent
\indent
\, $j$ \,\,\, times in the evaluation of \, $\phi(\f m{p_{i_u}}), 1 \le u \le j$, 
\newline
\indent
\indent
\indent
\indent
\indent
$\binom j2$ \, times in the evaluation of \, $\phi(\f m{p_{i_u}p_{i_v}}), 1 \le u < v \le j$,
\newline
\indent
\indent
\indent
\indent
\indent
$\binom j3$ \, times in the evaluation of \, $\phi(\f m{p_{i_u}p_{i_v}p_{i_w}}), 1 \le u < v < w \le j$,
\newline
\indent
\indent
\indent
\indent
\indent
\indent
\qquad $\vdots$
\newline
\indent
\indent
\indent
\indent
\indent
$\binom jj$ \, times in the evaluation of \, $\phi(\f m{p_{i_1}p_{i_2} \cdots p_{i_j}})$.

\noindent
Totally, they have been counted
\newline
\indent
\indent
\indent
\indent
\indent
$-j + \binom j2 - \binom j3 + \cdots + (-1)^j \binom jj = [(1-1)^j - 1] =-1$

\noindent
times.  Therefore, $\Phi_1(m, \phi)$ is indeed the number of periodic points of $f$ with least period $m$.  Similar argument applies to $\Phi_2$.  So, we obtain the following result:

\noindent
{\bf Theorem 1.}
{\it Let $S$ be a nonempty set and let $g$ be a map from $S$ into itself such that, for every positive integer $m$, the equation $g^m(x) = x$ (or $g^m(x) = -x$ respectively) has only finitely many distinct solutions.  Let $\phi(m)$ (or $\psi(m)$ respectively) denote the number of these solutions.  Then, for every positive integer $m$, the following hold:
\begin{itemize}
\item[(1)]
The number of periodic points of $g$ with least period $m$ is $\Phi_1(m, \phi)$.  Consequently, $\Phi_1(m, \phi) \equiv 0$ (mod $m$).

\item[(2)]
If $0 \in S$ and $g$ is odd, then the number of symmetric periodic points (i.e. periodic points whose orbits are symmetric with respect to the origin) of $g$ with least period $2m$ is $\Phi_2(m, \psi)$.  Consequently, $\Phi_2(m, \psi) \equiv 0$ (mod $2m$).
\end{itemize}}

Successful applications of the above theorem depend of course on a knowledge of the function $\phi$ or $\psi$.  For continuous maps from a compact interval into itself, the method of symbolic representations as introduced in 
{\bf{\cite{du1,du2,du3}}} is very powerful in enumerating the numbers (and hence generating the function $\phi$ or $\psi$) of the fixed points of all positive integral powers of the maps.  However, to get simple recursive formulas for the function 
$\phi$ or $\psi$, an appropriate map must be chosen.  The method of symbolic representations is simple, powerful, and easy to use.  Once you get the hang of it, the rest is only routine.  See {\bf{\cite{du1,du2,du3}}} for some examples regarding how this method works.  In the following, we present some new sequences which are found neither in {\bf{\cite{sl}}} nor in  
"superseeker\@research.att.com".  Proofs of these results can be followed from those of {\bf{\cite{du1,du2,du3}}}.

\noindent
{\bf Theorem 2.}
{\it For integers $n \ge 4$ and $1 < m < n-1$, let $f_{m,n}(x)$ be the continuous map from $[1, n]$ onto itself defined by:   $f_{m,n}(1) = m+1$, $f_{m,n}(2) = 1$, $f_{m,n}(m) = m-1$, $f_{m,n}(m+1) = m+2$, $f_{m,n}(n-1) = n$, $f_{m,n}(n) = m$, and $f_{m,n}(x)$ is linear on $[j, j+1]$ for every integer $j$ with $1 \le j \le n-1$.  Also let $f(x)$ be the continuous map from $[1, 4]$ onto itself defined by: $f(1)=f(3)=4$, $f(2)=1$, $f(4)=2$, and $f(x)$ is linear on $[1,2]$, $[2,3]$, and on 
$[3,4]$.  For integers $n \ge 3$, we also define sequences $< a_{n,k} >$ as follows: 
$$
a_{n,k} = 
\begin{cases}
2^{k+1} - 1, & \text{for} \,\,\, 1 \le k \le n-1, \\
3a_{n,k-1} - \sum_{i=2}^{n-1} a_{n,k-i}, & \text{for} \,\,\, n \le k. \\
\end{cases}
$$
\noindent
Then the following hold:
\begin{itemize}
\item[(a)]
For any positive integer $k$, $a_{3,k}$ is the number of distinct fixed points of the map $f^k(x)$ in $[1, 4]$, and for any positive integer $k$, any integers $n \ge 4$ and $1 < m < n-1$, the number of distinct fixed points of the map $f_{m,n}^k(x)$ in $[1, n]$ is $a_{n,k}$ which is clearly independent of $m$ for all $1 < m < n-1$.  Consequently, for any
integer $n \ge 3$, if $\phi_{a_n}(k) = a_{n,k}$ and $\Phi_1$ is defined as in Theorem 1, then $\Phi_1(k,\phi_{a_n}) \equiv 0$ (mod $k$) for all integers $k \ge 1$.

\item[(b)]
For every integer $n \ge 3$, the generating function $G_{a_n}(z)$ of the sequence $<a_{n,k}>$ is $G_{a_n}(z) =
(3z - \sum_{k=2}^{n-1} kz^k)/(1 - 3z + \sum_{k=2}^{n-1} z^k)$.
\end{itemize}}

\noindent
{\bf Theorem 3.}
{\it For every integer $n \ge 1$, let $g_n(x)$ be the continuous map from $[1, 2n+1]$ onto itself defined by: $g_n(1) = n+1$, $g_n(2) = 2n+1$, $g_n(n+1) = n+2$, $g_n(n+2) = n$, $g_n(2n+1) = 1$, and $g_n(x)$ is linear on $[j, j+1]$ for every integer $j$ with $1 \le j \le 2n$.  We also define sequences $<b_{n,k}>$ as follows:
$$
\begin{cases}
b_{n,2k-1} = 1, & \text{for} \,\,\, 1 \le k \le n, \\
b_{n,2k-1} = 2^{k-n-1}(2k-1)+1, & \text{for} \,\,\, n+1 \le k \le 2n, \\
b_{n,2k} = 2^{k+1}-1, & \text{for} \,\,\, 1 \le k \le 2n, \\
b_{n,k} = 3b_{n,k-2} - \sum_{i=2}^{2n}b_{n,k-2i}, & \text{for} \,\,\, k \ge 4n+1. \\
\end{cases}
$$
\noindent
Then, for any integers $k \ge 1$ and $n\ge 1$, $b_{n,k}$ is the number of distinct fixed points of the map $g_n^k(x)$ in $[1, 2n+1]$.  Consequently, if $\phi_{b_n}(k) = b_{n,k}$ and $\Phi_1$ is defined as in Theorem 1, then $\Phi_1(k,\phi_{b_n}) \equiv 0$ (mod $k$) for all integers $k \ge 1$.  Moreover, the generating function $G_{b_n}(z)$ of the sequence $<b_{n,k}>$ is $G_{b_n}(z) = (z + \sum_{k=2}^{2n} (-1)^kkz^k)/(1 - z - \sum_{k=2}^{2n} (-1)^kz^k)$.}

\noindent
{\bf Remark.}
In Theorem 3, when $n = 1$, the sequence $<b_{n,k}>$ becomes the Lucas sequence: 1,3,4,7,11, $\cdots$. 

\noindent
{\bf Theorem 4.}
{\it For integers $n \ge 2$, $2 \le j \le 2n+1$, and $2 \le m \le 2n+1$, let $h_{j,m,n}(x)$ be the continuous map from $[1, 2n+2]$ onto itself defined by: $h_{j,m,n}(1) = j$, $h_{j,m,n}(x) = 1$ for all even integers $x$ in $[2,2n]$, $h_{j,m,n}(x) = 2n+2$ for all odd integers $x$ in $[3, 2n+1]$, $h_{j,m,n}(2n+2) = m$, and $h_{j,m,n}(x)$ is linear on $[j, j+1]$ for every integer $j$ with $1 \le j \le 2n+1$.  We also define sequences $<c_{j,m,n,k}>$ as follows:
$$
c_{j,m,n,k} =
\begin{cases}
2n+1, & \text{for} \,\,\, k = 1, \\
(2n+1)^2 - 2[2n - (j-m)], & \text{for} \,\,\, k = 2, \\
(2n+1)^3 - 6n[2n+1 - (j-m)], & \text{for} \,\,\, k = 3, \\
(2n+1)c_{j,m,n,k-1}-[2n-(j-m)]c_{j,m,n,k-2}-(j-m)c_{j,m,n,k-3}, & \text{for} \,\,\, k \ge 4. \\
\end{cases}
$$}
\noindent
Then, for any integers $n \ge 2$, $2 \le j \le 2n+1$, $2 \le m \le 2n+1$, and $k \ge 1$, $c_{j,m,n,k}$ is the number of distinct fixed points of the map $h_{j,m,n}^k(x)$ in $[1, 2n+2]$.  Consequently, if $\phi_{c_{j,m,n}}(k) = c_{j,m,n,k}$ and $\Phi_1$ is defined as in Theorem 1, then $\Phi_1(k,\phi_{c_{j,m,n}}) \equiv 0$ (mod $k$) for all integers $k \ge  1$.  Moreover, the generating function $G_{c_{j,m,n}}(z)$ of the sequence $<c_{j,m,n,k}>$ is $G_{c_{j,m,n}}(z) = \{ \, (2n+1)z - 2[2n - (j-m)]z^2 - 3(j-m)z^3 \, \}/\{ \, 1 - (2n+1)z + [2n - (j-m)]z^2 + (j-m)z^3 \, \}$.

\noindent
{\bf Remarks.} 
(1) For fixed integers $n \ge 2, q, r$, and $s$, let $\phi(k)$ be the map on the set of all positive integers defined by:  $\phi(1) = 2n+1$, $\phi(2)= (2n+1)^2 - 2q$, $\phi(3) = (2n+1)^3 - 6r$ and $\phi(k) = (2n+1)\phi(k-1) - q\phi(k-2) - s\phi(k-3)$ for all integers $k \ge 4$.  Then Theorem 4 implies that, for some suitable choices of $q$, $r$, $s$, and a map $f$, $\phi(k)$ are the numbers of fixed points of $f^k(x)$ and hence, for $\Phi_1$ defined as in Theorem 1, $\Phi_1(k, \phi) \equiv 0$ (mod $k$) for all integers $k \ge 1$.  If we only consider $\phi(k)$ as a sequence of positive integers and disregard whether it represents the numbers of fixed points of all positive integral powers of some map, we can still ask if $\Phi_1(k,\phi) \equiv 0$ (mod $k$) for all integers $k \ge 1$.  Extensive computer experiments suggest that this seems to be the case for some other choices of $q$, $r$, and $s$.  Therefore, there should be a number-theoretic approach to this more general problem as does in Theorem 5 below.

(2) Note that, in Theorem 4 above, when $j = 2$ nd $m = 2n+1$, we actually have $c_{2,2n+1,n,k} = (2n-1)^k + 2$ which satisfies the difference equation $c_{2,2n+1,n,k+1} = (2n-1)c_{2,2n+1,n,k} - 4(n-1)$ for all positive integers $k$.

The following result concerning the linear recurrence of second-order can be obtained by counting the fixed points of all positive integral powers of maps similar to those considered in Theorem 4.  The number-theoretic approach can also be found in {\bf{\cite{sh,br}}}.  

\noindent
{\bf Theorem 5.}
{\it For integers $n \ge 2$ and $1-n \le m \le n$, let $<d_{m,n,k}>$ be the sequences defined by
$$
d_{m,n,k} =
\begin{cases}
n, &\qquad \text{for} \,\,\, k = 1, \\
n^2 + 2m, &\qquad \text{for} \,\,\, k = 2, \\
nd_{m,n,k-1} + md_{m,n,k-2}, &\qquad \text{for} \,\,\, k \ge 3. \\
\end{cases}
$$
\noindent
For any integers $n \ge 2$, $1-n \le m \le n$ and $k \ge 1$, if $\phi_{d_{m,n}}(k) = d_{m,n,k}$ and $\Phi_1$ is defined as in Theorem 1, then $\Phi_1(k, \phi_{d_{m,n}}) \equiv 0$ (mod $k$) for all integers $k \ge 1$.  Moreover, the generating
function $G_{d_{m,n}}(z)$ of the sequence $<d_{m,n,k}>$ is $G_{d_{m,n}}(z) = (nz + 2mz^2)/(1 - nz - mz^2)$.}

The following result is taken from {\bf [4}, Theorem 3{\bf ]}.  More similar examples can also be found in {\bf{\cite{du2}}}.

\noindent
{\bf Theorem 6.}
{\it For every integer $n \ge 2$, let $p_n(x)$ be the continuous odd map from $[-n,n]$ onto itself defined by $p_n(i) = i+1$ for every integer $i$ with $1 \le i \le n-1$, $p_n(n) = -1$, and $p_n(x)$ is linear on $[j, j+1]$ for every integer $j$ with $-n \le j \le n-1$.  We also define sequences $< s_{n,k} >$ as follows:
$$
s_{n,k} =
\begin{cases}
1, &\qquad \text{for} \,\,\, 1 \le k \le n-1, \\
2^{k-n}(2k) + 1, &\qquad \text{for} \,\,\, n \le k \le 2n-1, \\
3s_{n,k-1} - \sum_{i=2}^{2n-1} s_{n,k-i}, &\qquad \text{for} \,\,\, 2n \le k. \\
\end{cases}
$$
\noindent
Then, for any integers $n \ge 2$ and $k \ge 1$, $a_{2n,k}$ is the number of distinct fixed points of the map $p_n^k(x)$ in $[-n, n]$, where $a_{2n,k}$ is defined as in Theorem 2, and $s_{n,k}$ is the number of distinct solutions of the equation $p_n^k(x) = -x$ in $[-n,n]$.  Consequently, if $\psi_{s_n}(k) = s_{n,k}$ and $\Phi_2$ is defined as in Theorem 1, then $\Phi_2(k, \psi_{s_n}) \equiv 0$ (mod $2k$).  Moreover, the generating function $G_{s_n}(z)$ of $<s_{n,k}>$ is  $G_{s_n}(z) = [z-2z^2-z^3+ \sum_{k=5}^{n-1} (k-4)z^k + (3n-4)z^n - \sum_{k=n+1}^{2n-1} (2n-k)z^k]/(1- 3z + \sum_{k=2}^{2n-1} z^k)$.
(When $n = 2$, ignore $-2x^2$, and when $n = 3$, ignore $-x^3$).}

\noindent
{\bf Remark.} 
Numerical computations suggest that the maps $\psi_{s_n}$ in Theorem 6 also satisfy $\Phi_1(k, \psi_{s_n}) \equiv 0$ (mod $k$) for all integers $k \ge 1$.  However, our method cannot verify this.  There may be an algebraic-theoretic verification of it.

\noindent
{\bf ACKNOWLEDGMENTS}
The author wants to thank Professor Peter Jau-Shyong Shiue for his many invaluable suggestions and encouragements in writing this paper.

\end{document}